\documentclass[aps,reprint,twocolumn,pre,superscriptaddress]{revtex4-2}
\usepackage{graphicx}
\usepackage[usenames,dvipsnames]{xcolor}
\usepackage{amsfonts,amsbsy,color}
\usepackage{bm}
\usepackage{amsmath,amssymb}
\usepackage{mathtools}
\usepackage{tikz,tikz-cd}
\usepackage{jabbrv}

\usepackage[hidelinks]{hyperref}

\hypersetup{
  colorlinks   = true, %Colours links instead of ugly boxes
  urlcolor     = blue, %Colour for external hyperlinks
  linkcolor    = blue, %Colour of internal links
  citecolor   = blue %Colour of citations, could be ``red''
}

\newcommand{\bbE}{\mathbb{E}}

\newcommand{\bbR}{\mathbb{R}}

\newcommand{\bms}{\bm{s}}

\newcommand{\bmS}{\bm{S}}

\newcommand{\rmc}{\mathrm{c}}

\newcommand{\rmJ}{\mathrm{J}}

\newcommand{\sfJ}{\mathsf{J}}

\begin{document}
\title{Detailed balance and phase transition in Nash equilibria of mean-field Ising game}

\author{Ken Hiura}
\email{ken.hiura@yamaguchi-u.ac.jp}
\affiliation{Graduate School of Sciences and Technology for Innovation, Yamaguchi University, Yamaguchi 753-8512, Japan}

\author{Hajime Koike}
\email{koike.hajime.7i@kyoto-u.ac.jp}
\affiliation{Graduate School of Informatics, Kyoto University, Kyoto 606-8501, Japan}

\author{Kiyoshi Kanazawa}
\email{kiyoshi@scphys.kyoto-u.ac.jp}
\affiliation{Department of Physics, Graduate School of Science, Kyoto University, Kyoto 606-8502, Japan}

% \date{\today}

\begin{abstract}
    We study a finite-state mean-field game with entropic control cost. In the large population limit, we show that its stationary Nash equilibria satisfy the detailed balance condition and their distributions take the Gibbsian form with the value function playing the role of a Hamiltonian. We apply the entropic control cost to a ferromagnetic mean-field Ising game and investigate its stationary and non-stationary Nash equilibria. We find that the stationary Nash equilibria exhibit a second-order phase transition from disordered to ordered equilibria as temperature decreases. A center-manifold reduction yields non-stationary Nash equilibrium near the critical temperature and reveals the critical slowing down in the relaxation toward the stationary Nash equilibrium.
\end{abstract}

\maketitle

%Introduction

\section{Introduction}

Game theory provides a framework for describing systems of interacting decision makers whose actions are chosen strategically. A central concept in game theory is the Nash equilibrium, in which no player can reduce their own cost by changing their strategy unilaterally. When the number of players is very large, however, a direct analysis of the Nash equilibria becomes intractable. Mean-field game (MFG) theory \cite{LasryLions2007,HuangMalhameCaines2006} provides an effective description of such large population games by replacing the microscopic influence of the other players with a population distribution. Since its introduction, MFG has been intensively developed \cite{Carmona2018a, Carmona2018b} and applied to a wide range of problems, including economics and finance \cite{ChanSircar2017, Carmona2015}, epidemics \cite{LaguzetTurinici2015,AurellCarmonaDayanikliLauriere2022}, and biological phenomena \cite{GroverBakshiTheodorou2018, CarmonaGraves2020}, to cite a few.

Multiple Nash equilibria are common in mean-field games once the standard monotonicity conditions ensuring uniqueness are violated, and the selection problem has been extensively studied, e.g., \cite{BardiMarkus2019,CecchinPraFischerPelino2019,Delarue2019, CecchinDelarue2021}. A particularly interesting situation arises when the structure of multiple Nash equilibria changes qualitatively as a model parameter is varied. Such bifurcations provide a game-theoretic analogue of phase transitions in statistical physics. A pioneering work \cite{YinMehtaMeynShanbhag2011} and subsequent studies \cite{CarmonaCormierSoner2023, CarmonaCormierSoner2025} investigated transitions from incoherent to synchronized equilibria in game-theoretic counterparts of mean-field XY and Kuramoto models. More recently, a continuous bifurcation in a two-state MFG, closely analogous to the second-order transition of the Ising model, was analyzed by \cite{HoferSoner2026}.

In these works \cite{YinMehtaMeynShanbhag2011,CarmonaCormierSoner2023}, the MFGs are formulated in terms of controlled diffusion processes with quadratic control costs. In this setting, the best-response drift of each player against other players is proportional to the negative gradient of the value function. At a stationary Nash equilibrium, the controlled process satisfies detailed balance, and its stationary distribution takes a Gibbsian form \cite{LasryLions2006, CarmonaCormierSoner2023}, with the value function playing the role of a Hamiltonian in statistical mechanics. Although an analogous structure is known for single-agent optimal control of finite-state jump processes with relative-entropy costs \cite{TheodorouTodorov2012, ChetriteTouchette2015, Gueant2020, GaoLiuTse2023, HoriguchiKobayashi2025}, to our knowledge, it has not been explicitly discussed for stationary Nash equilibria in finite-state MFGs.

In this paper, we study a finite-state mean-field game with an entropic control cost, given by a relative entropy rate, penalizing deviations from the uncontrolled infinite-temperature process. We show that at a stationary MFG Nash equilibrium, the controlled jump process satisfies detailed balance and the stationary distribution is expressed as the Gibbs form with the value function in place of the Hamiltonian. To uncover the collective consequences of the MFG with the entropic cost, we propose a ferromagnetic Ising game with all-to-all interactions and the entropic cost. Each player controls the flip rate of a binary spin and pays two competing running costs: an energetic cost favoring ferromagnetic ordering and the entropic control cost. We reduce the problem of finding a Nash equilibrium to solving a two-dimensional dynamical system and analyze its local stability and bifurcation structure. We find that the stationary MFG Nash equilibria with nonzero magnetization $m_*(\beta)$ bifurcate continuously from the disordered equilibrium at a critical inverse temperature $\beta_{\rmc}$ with $m_{*}(\beta) \sim (\beta - \beta_{\rmc})^{\frac{1}{2}}$. The resulting exponent is the mean-field Ising value $\frac{1}{2}$, in contrast to the exponent $1$ found by Ref.~\cite{HoferSoner2026} in a two-state MFG with a quadratic cost. To study non-stationary Nash equilibria near the critical temperature, using the center-manifold reduction, we derive a reductive expression of the non-stationary Nash equilibria and study its critical dynamical properties. We show that the relaxation time $\tau (\beta)$ to the stationary Nash equilibrium exhibits the critical slowing down with $\tau (\beta) \sim |\beta - \beta_{\rmc}|^{-1}$.

The remainder of the paper is organized as follows. In Section \ref{sec:general}, we study a generic structure of finite-state mean-field games with entropic control costs. In Section \ref{sec:ising}, we introduce a mean-field Ising game and analyze the stationary MFG Nash equilibria and the associated phase transition and critical dynamics. Section \ref{sec:conclusion} concludes the paper. Technical derivations are collected in the Appendices.

\section{General structure}
\label{sec:general}

We demonstrate the general setup and properties of finite-player games with entropic control costs and their mean-field limits using finite-state jump processes. For diffusion processes, see comprehensive textbooks \cite{Carmona2018a,Carmona2018b}.

\subsection{Finite-player game}

We define a finite-state jump game with finitely many players. Consider $N$ players indexed by $i \in [N] \coloneqq \{1,\dots,N\}$. Each player has a spin random process. The state of the $i$-th player at time $t$ is given by $S_i(t)$, which takes values in a finite set $E$. We use $\bmS = (S_1,\dots, S_N)$ and $\bmS_{-i}$ to denote a collection of $N$ spins and the collection with the $i$-th component removed, respectively. We also introduce $(\bmS_{-i},S) = (S_1,\dots, S_{i-1},S,S_{i+1},\dots,S_N)$ to represent a concatenation of $\bmS_{-i}$ and $S$. 

The collection of spins $\bmS$ is a controlled Markov jump process. The only possible jumps are single spin state change events. The $i$-th player can act on the system by choosing the transition rate of the $i$-th spin state change. The controlled transition rate $\lambda_i(t,s,s^{\prime})$ at time $t$ associated with the $i$-th spin state change $s \to s^{\prime}$ is parameterized by the $i$-th player's control parameter $\alpha^i(t) :E^2 \backslash \Delta_E \to \bbR$ in the exponential tilting form,
\begin{align} \label{eq:lambda}
    \lambda_i(t,s,s^{\prime}) = k(s,s^{\prime}) e^{\alpha_{i}(t,s,s^{\prime})},
\end{align}
where $k$ is a transition rate in an uncontrolled process and $\Delta_E \coloneqq \{ (s,s) \in E^2 \}$ is the diagonal set of $E$. The admissible control is restricted to Markovian feedback controls, $\alpha_i(t) = \varphi_i (t,\bmS(t))$ for a deterministic function $\varphi_i$. Thus, the control parameter of the $i$-th player at time $t$, $\alpha_i(t)$, can depend on the current values of the states of the other players $\bmS_{-i}(t)$, as well as the current value of the state of the $i$-th player $S_i(t)$.

The objective of each player is to minimize their costs as much as possible. The cost function $C_i$ for the $i$-th player is defined by
\begin{align} \label{eq:cost}
    C_i(t,\bms,\bm{\alpha}) \coloneqq \bbE \left[ \int_t^{\infty} du \ e^{-r(u-t)} c_i(\bmS(u),\alpha_i(u)) \right],
\end{align}
with the initial condition $\bmS (t) = \bms$. The discount factor $e^{-r(u-t)}$ ($r > 0$) is introduced to ensure the convergence of the integral. The running cost function $c$ is composed of two running costs $c_{i1}$ and $c_2$:
\begin{align}
    c_i(\bms,\alpha) = c_{i1}(\bms) + 2 \beta^{-1} c_{2}(s_i,\alpha).
\end{align}
The first cost $c_{i1}$ quantifies the ``energetic'' running cost depending on only the spin configuration. The second cost $c_2$ quantifies the ``entropic'' cost associated with deviating from a reference process due to control. The proportionality constant $\beta^{-1} > 0$ in front of the entropic cost plays the role of the temperature of the game and controls the competition between energetic and entropic costs: the entropic cost is relatively dominant at high temperatures $\beta \ll 1$, while the energetic cost is dominant at low temperatures $\beta \gg 1$. In this paper, the energetic cost may vary depending on the model under consideration, whereas the entropic cost is chosen to satisfy the Gibbsian structure and the detailed balance condition at stationary Nash equilibria independently of the model. We provide a detailed discussion of the choice of $c_2$ in the next section.

\subsection{Entropic control cost}

We choose the control cost $c_2$ as the instantaneous relative entropy rate associated with the deviation of an individual player's controlled jump process from an uncontrolled reference process. Let $p_{(\bm{\alpha}_{-i},\alpha)}(t+h,\bms^{\prime}|t,\bms)$ and $p_{(\bm{\alpha}_{-i},0)}(t+h,\bms^{\prime}|t,\bms)$ denote the short-time transition probability density from $\bms$ at time $t$ to $\bms^{\prime}$ at time $t+h$ of the controlled processes when the $i$-th player uses the control $\alpha$ and the reference control $\alpha = 0$ at time $t$, respectively. The control cost $c_2$ is chosen to be the instantaneous relative entropy rate of the controlled process with respect to the reference process conditioned on the current state being $\bms$ with $s_i = s$,
\begin{align}
    c_2(s,\alpha) &\coloneqq \lim_{h \downarrow 0} \frac{D(p_{(\bm{\alpha}_{-i},\alpha)}(t+h,\cdot|t,\bms)\|p_{(\bm{\alpha}_{-i},0)}(t+h,\cdot|t,\bms))}{h},
\end{align}
where the right-hand side is actually independent of $i$ and $\bms_{-i}$. In the remainder of the paper, we take $k(s,s^{\prime}) = 1$ for all $(s,s^{\prime}) \in E^2 \backslash \Delta_E$. The reference process has the uniform invariant distribution and may be regarded as an infinite-temperature process. Consequently, we take
\begin{align} \label{eq:entcost}
    c_{2}(s,\alpha) \coloneqq \sum_{s^{\prime} \in E \backslash \{s\}} \left( e^{\alpha(s,s^{\prime})} \alpha (s,s^{\prime}) - e^{\alpha (s,s^{\prime})} + 1 \right).
\end{align}
The derivation of Eq.~\eqref{eq:entcost} is given in Appendix \ref{app:entcost}. Notice that this cost is nonnegative and vanishes if and only if $\alpha(s,s^{\prime}) = 0$ for every transition, as expected.

This choice is a counterpart in jump process of the quadratic cost for controlled diffusion processes. Indeed, the quadratic control cost commonly used for the controlled diffusion processes has a natural interpretation in terms of the relative entropy. Consider an $\bbR^n$-valued controlled diffusion process
\begin{align} \label{eq:SDE}
    dS_i(t) = \alpha_i(t) dt + dW_i(t),
\end{align}
where $\alpha_i$ is the controlled drift and $W_i$ is an $n$-dimensional Wiener process. A standard choice of the running control cost is the quadratic cost,
\begin{align} \label{eq:SDEentcost}
    c_2(s,\alpha) = \frac{1}{2} \alpha^2,
\end{align}
which penalizes deviations of the controlled diffusion process from the uncontrolled free diffusion. On the other hand, the instantaneous relative entropy rate is given by
\begin{align}
    \lim_{h \downarrow 0} \frac{D(p_{(\bm{\alpha}_{-i},\alpha)}(t+h,\cdot|t,\bms)\|p_{(\bm{\alpha}_{-i},0)}(t+h,\cdot|t,\bms))}{h} = \frac{1}{2} \alpha^2.
\end{align}
Thus, the quadratic running cost \eqref{eq:SDEentcost} is precisely the instantaneous relative entropy rate of the controlled diffusion with respect to the uncontrolled diffusion. The uncontrolled free diffusion process $dS_i = dW_i$ corresponds to the uncontrolled jump process with $k = 1$ in the above formulation.

An alternative choice of the control cost in finite-state jump games is a quadratic cost $c_2 = \frac{1}{2}\alpha^2$ with the additive parametrization $\lambda = 1 + \alpha$, e.g., \cite{HoferSoner2026}. Such a construction, however, does not preserve the relative-entropy interpretation of the quadratic cost in the diffusion processes. In this paper, for finite-state jump games, we instead use the instantaneous relative entropy rate as the control cost and study its consequences. We note that for weak controls $|\alpha| \ll 1$, $\lambda \simeq 1 + \alpha$ in Eq.~\eqref{eq:lambda} and
\begin{align}
    c_2(s,\alpha) \simeq \sum_{s^{\prime} \in E \backslash \{s\}} \frac{1}{2} \alpha(s,s^{\prime})^2.
\end{align}
The entropic cost given by Eq.~\eqref{eq:entcost} therefore reduces to the usual quadratic cost for weak controls, while retaining the relative-entropy interpretation for finite controls.

\subsection{Nash equilibria}

A fundamental notion of optimality in game theory is the Nash equilibrium. In this paper, we investigate Markov perfect Nash equilibria, which are the Markovian analogue of subgame perfect equilibria. A collection of Markov feedback controls $\bm{\alpha}^* = (\alpha_1^*,\dots,\alpha_N^*)$ is said to be a Markov perfect Nash equilibrium, or simply, Markov Nash equilibrium of the game if, for every $(t,\bms)$, every $i \in [N]$ and every Markov feedback control of the $i$-th player $\alpha_i$,
\begin{align}
    C_i(t,\bms,\bm{\alpha}^*) \leq C_i(t,\bms,(\bm{\alpha}_{-i}^*,\alpha_i)).
\end{align}
Intuitively, this is an equilibrium state in which each player is unable to reduce their own cost by changing their strategy alone.

A Markov Nash equilibrium is characterized by a fixed point of the best response function. Let us fix strategies $\bm{\alpha}_{-i}$ and consider the optimal control problem of the $i$-th player. For $(t,\bms)$, we define $v_i^N$, called a value function, as
\begin{align}
    v_i^N(t,\bms ; \bm{\alpha}_{-i}) \coloneqq \inf_{\varphi} C_i(t,\bms,(\bm{\alpha}_{-i},\varphi)),
\end{align}
where the minimization is taken over the class of Markov feedback controls. This is the cost when the $i$-th player responds in the best way to the other players' strategies. Owing to the particular relative-entropy form of the control cost, the minimization with respect to the flip control can be performed explicitly. The minimizing control is
\begin{align} \label{eq:finiteNash}
    &\varphi_i^{N,*}(t,\bms ; \bm{\alpha}_{-i})(s_i,s^{\prime})
    \notag \\
    &= - \frac{\beta}{2}\left[v_i^N(t,(\bms_{-i},s^{\prime}) ; \bm{\alpha}_{-i})-v_i^N(t,\bms ; \bm{\alpha}_{-i})\right].
\end{align}
The dynamic programming principle gives a Hamilton--Jacobi--Bellman (HJB) equation for the value function. By definition, $\bm{\alpha}^*(t) = \bm{\varphi}^*(t,\bmS(t))$ is a Markov Nash equilibrium if and only if $\varphi_i^{N,*}(t,\bms ; \bm{\alpha}_{-i}^*) = \varphi_i^*(t,\bms)$ for every $i \in [N]$. Imposing this fixed-point condition that every player simultaneously uses their best response yields the coupled HJB equation
\begin{align} \label{eq:cHJB}
    0={}&\partial_t v_i^N-r v_i^N + c_{i1}(\bms)
    \notag\\
    &+2\beta^{-1}\sum_{s^{\prime}\in E\backslash\{s_i\}}
    \left[1-e^{-\frac{\beta}{2}\left(v_i^N(t,(\bms_{-i},s^{\prime}))-v_i^N(t,\bms)\right)}\right]
    \notag\\
    &+\sum_{j\in[N]\backslash\{i\}}\sum_{s^{\prime}\in E\backslash\{s_j\}}
    e^{-\frac{\beta}{2}\left(v_j^N(t,(\bms_{-j},s^{\prime}))-v_j^N(t,\bms)\right)}
    \notag\\
    &\qquad\qquad\times\left[v_i^N(t,(\bms_{-j},s^{\prime}))-v_i^N(t,\bms)\right],
\end{align}
for $i\in[N]$. Therefore, once the value functions are obtained from the coupled HJB equation, Eq.~\eqref{eq:finiteNash} gives the Markov Nash equilibrium. A derivation of Eqs.~\eqref{eq:finiteNash} and \eqref{eq:cHJB} is given in Appendix \ref{app:cHJB}.

\subsection{Mean-field game}

The coupled HJB equation \eqref{eq:cHJB} is difficult to solve for a generic energetic cost $c_{i1}$ and for a finite $N$. We consider the large population limit $N \uparrow \infty$ for a particular class of the energetic cost $c_{i1}$. We now suppose that the influence on the $i$-th spin from the other spins is effectively described by the mean-field. Namely, $c_{i1}$ is supposed to have the following form:
\begin{align} \label{eq:MFcost}
    c_{i1}(\bms) = \bar{c}(s_i, \rho^N_{\bms,i}) + O(N^{-1}),
\end{align} 
where $\rho^N_{\bms,i}$ is the empirical density function of all players except $i$ defined by
\begin{align}
    \rho^N_{\bms,i}(s)
    \coloneqq \frac{1}{N-1}\sum_{j\in[N]\backslash\{i\}}\delta_{s,s_j}.
\end{align}
This form of the energetic running cost makes it possible to take the mean-field limit.

The mean-field limit of the coupled HJB equation is reduced to the so-called MFG system. Let $\rho(t,\cdot)$ denote a probability density function on $E$ at time $t$. The MFG system is described by an HJB equation for the representative player's value function $v(t,s)$ and a forward master equation for $\rho(t,s)$:
\begin{align} \label{eq:MFHJB}
    0={}&\partial_t v(t,s)-r v(t,s) + \bar{c}(s,\rho(t))
    \notag\\
    &+2\beta^{-1}\sum_{s^{\prime}\in E\backslash\{s\}}
    \left[1-e^{-\frac{\beta}{2}\left(v(t,s^{\prime})-v(t,s)\right)}\right],
\end{align}
and
\begin{align} \label{eq:MFmaster}
    0={}&\partial_t\rho(t,s)
    -\sum_{s^{\prime}\in E\backslash\{s\}}
    \left[\rho(t,s^{\prime})e^{-\frac{\beta}{2}\left(v(t,s)-v(t,s^{\prime})\right)}\right.
    \notag\\
    &\left. \phantom{\partial_t \rho(t,s)} -\rho(t,s)e^{-\frac{\beta}{2}\left(v(t,s^{\prime})-v(t,s)\right)}\right]
\end{align} 
with the initial condition $\rho (0,s) = \rho_0(s)$, together with the boundedness condition for $v$ discussed below. 

If the initial spins $\bmS (0)$ are sampled independently according to $\rho_0$, the finite-player Nash equilibrium \eqref{eq:finiteNash} along typical trajectories can be approximated as
\begin{align} \label{eq:approxNash}
    &\varphi_i^{N,*}(t,\bmS(t))(S_i(t),s^{\prime})
    \notag \\
    &\approx -\frac{\beta}{2}\left[v(t,s^{\prime})-v(t,S_i(t))\right].
\end{align}
Thus, by solving the MFG system, \eqref{eq:MFHJB} and \eqref{eq:MFmaster}, under the initial condition $\rho (0,\cdot) = \rho_0$ and the boundedness condition for $v$, we find an approximate Markov Nash equilibrium in the large population limit. A heuristic derivation of Eqs.~\eqref{eq:MFHJB}, \eqref{eq:MFmaster}, and \eqref{eq:approxNash} is given in Appendix \ref{app:MFG}.

The MFG system can be formulated as the optimal control theory for the representative player with the self-consistency condition. The objective of the representative player is to minimize the cost defined by
\begin{align}
    C(t,s,\alpha,\rho) &= \bbE \left[ \int_t^{\infty} du \ e^{-r(u-t)} \left( \bar{c}(S(u),\rho(u)) \right. \right.
    \notag \\
    &\phantom{\bbE [ \int_t^{\infty} du]}\left. \left. + 2 \beta^{-1} c_2(S(u),\alpha(u)) \right) \right],
\end{align}
under the initial condition $S(t) = s$ and a given flow of probability density function $\{ \rho (u) \}_{u \geq 0}$. We define the value function $v(t,s,\rho) = \inf_{\alpha} C(t,s,\alpha,\rho)$ for this problem. Then, the dynamic programming principle leads to the HJB equation \eqref{eq:MFHJB}. The optimal feedback control is given by
\begin{align}
    \alpha^* (t,s,s^{\prime},\rho) = - \frac{\beta}{2} \left[ v(t,s^{\prime},\rho) - v(t,s,\rho) \right].
\end{align}
The MFG system is obtained by imposing the self-consistency condition that the prescribed population flow coincides with the law of the spin variable of the representative player under the optimal feedback control. Namely, $\rho$ is required to satisfy
\begin{align}
    \partial_t \rho = \sum_{s^{\prime} \in E \backslash \{s\}} \left[ \rho (t,s^{\prime}) e^{\alpha^* (t,s^{\prime},s,\rho)} - \rho (t,s) e^{\alpha^*(t,s,s^{\prime},\rho)} \right]
\end{align}
with $\rho (0,s) = \rho_0(s)$. Then, $(v(\rho),\rho)$ solves Eqs.~\eqref{eq:MFHJB} and \eqref{eq:MFmaster}.

The optimal control formulation provides the natural boundary condition for the HJB equation. Appendix \ref{app:bound} shows that if $\bar{c}$ is a continuous function of $\rho$, the value function is uniformly bounded in time,
\begin{align}
    \sup_{t \geq 0} |v(t,s)| < \infty.
\end{align}
Therefore, we impose this boundedness condition on the value function when we solve the differential equation \eqref{eq:MFHJB} and \eqref{eq:MFmaster}.

\subsection{Gibbsian structure and detailed balance}

A solution of Eqs.~\eqref{eq:MFHJB} and \eqref{eq:MFmaster} under the boundary condition is called an MFG Nash equilibrium. In this section, we study particularly stationary solutions of Eqs.~\eqref{eq:MFHJB} and \eqref{eq:MFmaster}, called stationary MFG Nash equilibria.

A useful consequence of the entropic control cost is that the transition rates and the probability distribution in a stationary Nash equilibrium acquire a Gibbsian structure and detailed balance condition. Indeed, for a time-independent value function $v$, the transition rate at the Nash equilibrium is given by
\begin{align}
    \lambda^* (s,s^{\prime}) = e^{- \frac{\beta}{2} \left[ v(s^{\prime}) - v(s) \right]}.
\end{align}
This transition rate satisfies the detailed balance condition with respect to the Gibbs distribution with the value function $v$ in place of the Hamiltonian,
\begin{align} \label{eq:gibbs}
    \rho (s) = \frac{e^{- \beta v(s)}}{Z}, \quad Z \coloneqq \sum_{s \in E} e^{- \beta v(s)}.
\end{align}
Namely,
\begin{align}
    \rho (s^{\prime}) \lambda^* (s^{\prime},s) - \rho (s) \lambda^* (s,s^{\prime}) = 0
\end{align}
for every transition $(s,s^{\prime}) \in E^2 \backslash \Delta_E$. This Gibbsian structure provides a natural interpretation of $\beta^{-1}$ as an effective temperature, with the stationary value function playing the role of an effective energy.

This Gibbsian structure and the detailed balance condition appear in the diffusion mean-field games with quadratic control cost \eqref{eq:SDEentcost}; see, for example, \cite{CarmonaCormierSoner2023, LasryLions2006}. In the present finite-state jump setting, the same structure follows from the relative-entropy control cost \eqref{eq:entcost}. 

\section{Model: mean-field Ising game}
\label{sec:ising}

We investigate a mean-field Ising game with the entropic control cost. This is a game-theoretic version of the ferromagnetic mean-field Ising model and gives a prototypical example of a second-order phase transition of Nash equilibria in mean-field games. 

\subsection{Definition}

We define a ferromagnetic Ising game with all-to-all interactions. The state space is $E = \{ -1, +1 \}$. The energetic cost $c_{i1}$ is chosen to be
\begin{align} \label{eq:penecost}
    c_{i1}(\bms) \coloneqq - \frac{J}{2N} \left( \sum_{j \in [N] \backslash \{i\}} s_j \right) s_i,
\end{align}
Here, $J >0$ is a positive coupling constant. This cost corresponds to two-body interaction energy per spin of the ferromagnetic Ising model. In terms of the empirical magnetization
\begin{align}
    m^N_{\bms,i}
    =\sum_{s\in E}s\rho^N_{\bms,i}(s),
\end{align}
the energetic cost \eqref{eq:penecost} satisfies the mean-field asymptotic condition \eqref{eq:MFcost}:
\begin{align}
    c_{i1}(\bms)
    = - \frac{J}{2} m^N_{\bms,i} s_i + O(N^{-1}).
\end{align}

The mean-field Ising game may be interpreted as a model of opinion formation in a community. The $i$-th agent holds either a positive ($s_i = +1$) or negative ($s_i = -1$) opinion on a given issue. The ferromagnetic interaction \eqref{eq:penecost} represents the pressure to conform to the majority opinion, while the entropic control cost \eqref{eq:entcost} represents the effort required to deliberately regulate one's opinion changes. The temperature $\beta^{-1}$ controls the competition between these two effects. At high temperatures, agents tend to give up actively regulating their opinions to conform to the majority opinions and instead allow their opinions to fluctuate randomly. As a result, the community may settle into an equilibrium state in which the population is evenly divided between two opinions. At low temperatures, by contrast, agents are more willing to bear the effort of regulating their opinions in order to conform to the majority opinion, which may lead to an equilibrium state in which one of the two opinions becomes predominant. In what follows, we confirm this picture by analyzing the stationary Nash equilibria and the relaxational behavior toward them.

\subsection{MFG Nash equilibria}

Eqs.~\eqref{eq:MFHJB} and \eqref{eq:MFmaster} with this choice determine MFG Nash equilibria of the mean-field Ising game. Using
\begin{align}
    m(t) &= \rho (t,+1) - \rho (t,-1),
    \\
    \Delta (t) &= \frac{\beta}{2} \left( v(t,+1) - v(t,-1) \right),
\end{align}
we write Eqs.~\eqref{eq:MFHJB} and \eqref{eq:MFmaster} equivalently as a two-dimensional differential equation
\begin{align} \label{eq:mdot}
    \dot{m} &= - 2 \left( \sinh \Delta + m \cosh \Delta \right),
    \\ \label{eq:Deltadot}
    \dot{\Delta} &= r \Delta + \frac{\beta J}{2} m + 2 \sinh \Delta
\end{align}
with the initial condition $m (0) = m_0$ and the boundedness condition
\begin{align}
    \sup_{t \geq 0} |\Delta (t)| < \infty,
\end{align}
which follows from the boundedness condition for $v$. For a solution $(m,\Delta)$ of Eqs.~\eqref{eq:mdot} and \eqref{eq:Deltadot}, the value function and the probability density function are obtained in terms of $(m,\Delta)$ as
\begin{align}
    v (t,\pm 1) = \bar{v}(t) \pm \beta^{-1} \Delta (t),
\end{align}
and
\begin{align}
    \rho (t,\pm 1) = \frac{1 \pm m(t)}{2},
\end{align}
respectively. Here,
\begin{align}
    \bar{v}(t) = 2 \beta^{-1} \int_t^{\infty} du \ e^{-r(u-t)} (1- \cosh \Delta (u)).
\end{align}
The corresponding feedback control is
\begin{align}
    \alpha^*(t,+1, -1) = \Delta (t), \quad \alpha^*(t,-1,+1) = - \Delta (t).
\end{align}

\subsection{Second-order phase transition in stationary MFG Nash equilibria}

We first study the stationary MFG Nash equilibria. By eliminating $\Delta$ from Eqs.~\eqref{eq:mdot} and \eqref{eq:Deltadot} at a fixed point, we obtain a single self-consistency equation
\begin{align} \label{eq:Fm}
    F_{\beta} (m) \coloneqq \frac{r}{2} \ln \frac{1+m}{1-m} - \frac{\beta J}{2} m + \frac{2m}{\sqrt{1-m^2}} = 0.
\end{align}
Thus, each root $m \in (-1,1)$ of Eq.~\eqref{eq:Fm} determines a stationary MFG Nash equilibrium.

We now analyze how the solutions of Eq.~\eqref{eq:Fm} change with temperature. The disordered equilibrium $m=0$ is always a solution at all temperatures. The expansion of $F_{\beta}(m)$ around $m=0$ gives
\begin{align} \label{eq:Fexpandp2}
    F_{\beta}(m)
    =\frac{J}{2}(\beta_{\rm c}-\beta)m
    +\left(\frac{r}{3}+1\right)m^3+O(m^5),
\end{align}
where the critical inverse temperature $\beta_{\rmc}$ is defined by
\begin{align} \label{eq:betac}
    \beta_{\rm c} \coloneqq \frac{2(r+2)}{J}.
\end{align}
This local analysis suggests the existence of a second-order phase transition at $\beta = \beta_{\rmc}$. In Appendix \ref{app:beta}, we show that for $\beta<\beta_{\rm c}$, $m=0$ is the only stationary solution: at $\beta = \beta_{\rmc}$, a pitchfork bifurcation occurs; and for $\beta > \beta_{\rm c}$, two nonzero solutions $\pm m_*(\beta)$ emerge. Close to the critical point,
\begin{align} \label{eq:msp}
    m_*(\beta)
    \simeq
    \sqrt{\frac{J(\beta-\beta_{\rm c})}{2\left(\frac{r}{3}+1\right)}}
    \sim(\beta-\beta_{\rm c})^{1/2}.
\end{align}
The magnetization therefore vanishes continuously at the critical point with the mean-field exponent $\frac{1}{2}$, characteristic of the mean-field Ising universality class. Thus, the second-order phase transition emerges as a result of the strategic optimization of the decision makers, in contrast to statistical mechanics, where one starts from a prescribed probability distribution over microscopic configurations. 

A similar continuous transition in a two-state MFG was found in \cite{HoferSoner2026}. In our model, the energetic cost $c_1$ is in fact equivalent to the running cost considered in \cite{HoferSoner2026} up to an additive constant, whereas our relative-entropy control cost differs from the quadratic control cost used there. Nevertheless, on the one hand, the critical point coincides when the reference jump rate in \cite{HoferSoner2026} is set to $\sigma^2 = 1$. This agreement is not accidental. In our model, the controlled transition rate is parameterized as $\lambda = e^{\alpha}$, and hence $\lambda \simeq 1 + \alpha$ for weak controls. After the corresponding identification of variables and parameters, the two MFG systems linearized around the disordered equilibrium $m=0$ become identical and therefore yield the same linear instability threshold, $\beta_{\rmc}J = 2(r+2)$. On the other hand, in \cite{HoferSoner2026}, the magnetization scales linearly with the distance from the critical point, corresponding to a critical exponent $1$, while our model exhibits the mean-field Ising scaling with exponent $\frac{1}{2}$. This difference originates from the nonlinear structure of the MFG systems induced by the different control costs.

\begin{figure}[t]
 \centering
 \includegraphics[width=0.95\linewidth]{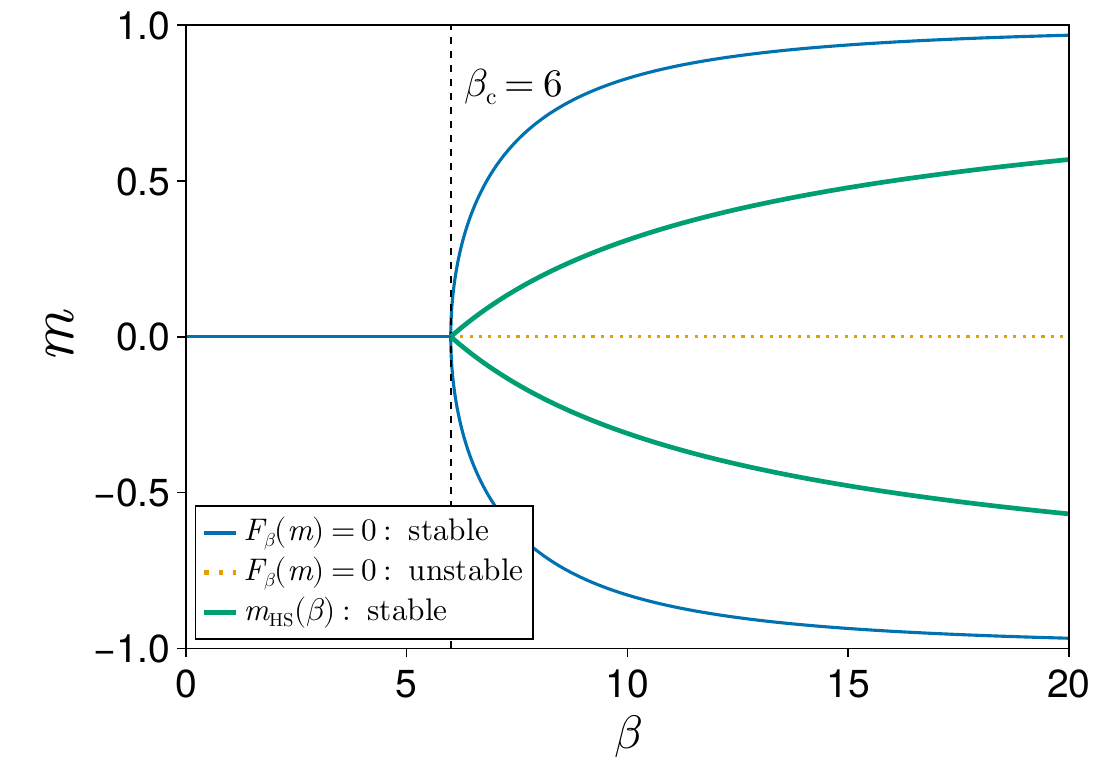}
 \caption{Stationary magnetization as a function of $\beta$ with $r=J=1$. Blue solid branches denote linearly stable solutions of $F_{\beta}(m) = 0$, while yellow dashed branches denote unstable solutions. The disordered equilibrium loses stability at $\beta_{\rmc}=6$, where the ordered branches bifurcate continuously. Green solid curves show the linearly stable stationary solutions $\pm m_{\mathrm{HS}}(\beta)$ of the model in \cite{HoferSoner2026} with reference jump rate $\sigma^2 = 1$.}
 \label{fig:p2}
\end{figure}

Figure \ref{fig:p2} shows the bifurcation diagram for $r=J=1$ obtained from numerically solving the algebraic equation \eqref{eq:Fm}. The stability classification follows from Sec. \ref{subsec:stability}. For comparison, the stable stationary solutions of the model in \cite{HoferSoner2026} with $\sigma^2 = 1$ are also shown.

\subsection{Linear stability}
\label{subsec:stability}

We analyze linear stability of stationary MFG Nash equilibria, taking into account that the initial magnetization is provided, whereas the initial value of $\Delta$ is chosen to satisfy the boundedness condition. As shown below, the linearization of Eqs.~\eqref{eq:mdot} and \eqref{eq:Deltadot} always has an unstable direction. Hence, stability of a stationary MFG Nash equilibrium cannot be characterized by ordinary Lyapunov stability under forward time evolution. Instead, a stationary MFG Nash equilibrium $(m_{*},\Delta_{*})$ is called linearly stable if, for every initial magnetization perturbation $\delta m(0)$, we can adjust the initial difference $\delta \Delta(0)$ such that the solution of the linearized MFG system for $(\delta m (t), \delta \Delta (t))$ decays to zero as $t \to \infty$. In other words, a one-dimensional stable eigenspace exists and has a nonzero projection onto the $m$-axis.

At a stationary equilibrium $(m_*,\Delta_*)$, we have $m_{*}= - \tanh \Delta_{*}$. The Jacobian matrix at $(m_{*},\Delta_{*})$ of Eqs.~\eqref{eq:mdot} and \eqref{eq:Deltadot} is
\begin{align}
    \mathsf{J} = 
    \begin{bmatrix}
        \displaystyle - \frac{2}{\sqrt{1-m_*^2}} & \displaystyle -2 \sqrt{1 - m_{*}^2} \\
        \displaystyle \frac{\beta J}{2}  & \displaystyle r + \frac{2}{\sqrt{1-m_{*}^2}}
    \end{bmatrix}.
\end{align}
Its trace and determinant are $\mathrm{tr} \mathsf{J} = r > 0$ and $\det \mathsf{J} = - 2 \sqrt{1-m_{*}^2} F_{\beta}^{\prime}(m_{*})$. If $F_{\beta}^{\prime}(m_{*}) > 0$, then the Jacobian has exactly one negative and one positive eigenvalue. Since $\sfJ_{12} \neq 0$, the stable eigenvector has a nonzero $m$-component. If $F_{\beta}^{\prime}(m_{*}) < 0$, both eigenvalues have positive real part and no nontrivial perturbation decays forward in time. Therefore, a fixed point $(m_{*},\Delta_{*})$ is linearly stable if and only if $F_{\beta}^{\prime}(m_{*}) > 0$. At the critical point, $F_{\beta}^{\prime}(m_{*}) = 0$ and the eigenvalues are $0$ and $r$.

The disordered equilibrium is stable for $\beta < \beta_{\rmc}$ and unstable for $\beta > \beta_{\rmc}$ since
\begin{align}
    F_{\beta}^{\prime}(0) = \frac{J}{2} (\beta_{\rmc} - \beta).
\end{align}
The two nonzero branches born at $\beta = \beta_{\rmc}$ satisfy
\begin{align} \label{eq:Fprime}
    F_{\beta}^{\prime}(m_{*}) = \frac{r}{2} \left( \frac{2}{1-m_{*}^2} - \frac{1}{m_{*}} \int_{-m_{*}}^{m_{*}} \frac{dx}{1+x} \right) + \frac{2m_{*}^2}{(1-m_{*}^2)^{\frac{3}{2}}}.
\end{align}
The strict convexity of $x \mapsto (1+x)^{-1}$ and the Hermite-Hadamard inequality imply that the first term is positive. Therefore, $\pm m_{*}(\beta)$ are linearly stable equilibria.

\subsection{Non-stationary Nash equilibrium near critical point}

The relaxational behavior to a linearly stable stationary MFG Nash equilibrium is described by a non-stationary Nash equilibrium solving Eqs.~\eqref{eq:mdot} and \eqref{eq:Deltadot}. Close to the critical temperature, we can compute perturbatively a non-stationary MFG Nash equilibrium near $(m,\Delta) = (0,0)$ using the center-manifold reduction.

At the critical point, the Jacobian matrix at $(m,\Delta) = (0,0)$ of Eqs.~\eqref{eq:mdot} and \eqref{eq:Deltadot} is
\begin{align}
    \sfJ_{\rmc} = 
    \begin{bmatrix}
        -2 & -2 \\
        r+2 & r+2
    \end{bmatrix}.
\end{align}
Let
\begin{align}
    \mathbf{u}_{+} =
    \begin{bmatrix}
        2 \\
        -(r+2)
    \end{bmatrix},
    \quad
    \mathbf{u}_{-} =
    \begin{bmatrix}
        1 \\
        -1
    \end{bmatrix}
\end{align}
be the right eigenvectors along the unstable and critical directions, respectively. Namely, $\sfJ_{\rmc} \mathbf{u}_{+} = r \mathbf{u}_{+}$ and $\sfJ_{\rmc} \mathbf{u}_{-} = 0$. The corresponding left eigenvectors are
\begin{align}
    \mathbf{v}_{+} = - \frac{1}{r}[1,1], \quad \mathbf{v}_{-} = \frac{1}{r} [r+2, 2].
\end{align}
We now express a non-stationary MFG Nash equilibrium as
\begin{align} \label{eq:cmfd}
    \mathbf{x} =
    \begin{bmatrix}
        m \\
        \Delta
    \end{bmatrix}
    = A \mathbf{u}_{-} + \rho (A) \mathbf{u}_{+}
\end{align}
and suppose
\begin{align} \label{eq:OmegaA}
    \dot{A} = \Omega (A).
\end{align}
Here, $\rho$ and $\Omega$ are determined perturbatively from the differential equation, \eqref{eq:mdot} and \eqref{eq:Deltadot}. Eqs.~\eqref{eq:mdot} and \eqref{eq:Deltadot} are expanded around $m=0$ as
\begin{align} \label{eq:ODE}
    \dot{\mathbf{x}} = \left( \rmJ_{\rmc} + \frac{\epsilon}{2} 
    \begin{bmatrix}
        0 & 0 \\
        1 & 0
    \end{bmatrix}
    \right)
    \mathbf{x}
    + \mathbf{N}(\mathbf{x}),
\end{align}
where $\epsilon \coloneqq (\beta - \beta_{\rmc}) J$ and
\begin{align}
    \mathbf{N}(\mathbf{x}) = 
    \begin{bmatrix}
        - \frac{\Delta^3}{3} - m \Delta^2 + O(\Delta^5, m \Delta^4) \\
        \frac{\Delta^3}{3} + O(\Delta^5)
    \end{bmatrix}.
\end{align}
By substituting Eq.~\eqref{eq:cmfd} into Eq.~\eqref{eq:ODE}, we obtain
\begin{align} \label{eq:invariance}
    \sfJ_{\rmc} \rho \mathbf{u}_{+} = \Omega \mathbf{u}_{-} + \Omega \rho^{\prime} \mathbf{u}_{+} - \frac{\epsilon}{2}
    \begin{bmatrix}
        0 & 0 \\
        1 & 0
    \end{bmatrix}
    \mathbf{x} - \mathbf{N}(\mathbf{x}).
\end{align}
The solvability condition $\langle \mathbf{v}_{-}, \sfJ_{\rmc} \rho \mathbf{u}_{+} \rangle = 0$ implies
\begin{align}
    \Omega &=  \left\langle \mathbf{v}_{-}, \frac{\epsilon}{2}
    \begin{bmatrix}
        0 & 0 \\
        1 & 0
    \end{bmatrix}
    \mathbf{x} + \mathbf{N}(\mathbf{x}) \right\rangle
    \notag \\ \label{eq:pOmega}
    &= \frac{\epsilon}{r} A - \frac{2(r+3)}{3r} A^3 + O(A^5).
\end{align}
By taking the inner product of Eq.~\eqref{eq:invariance} with $\mathbf{v}_{+}$ and assuming the pitchfork scaling $A = O(\epsilon^{\frac{1}{2}})$, we have
\begin{align}
    \rho = \frac{\epsilon}{2 r^2} A - \frac{A^3}{r^2} + O(A^5).
\end{align}
The critical properties near the critical point can be read from the reduced normal form \eqref{eq:pOmega} of the non-stationary Nash equilibrium. Above the critical temperature $\beta \lesssim \beta_{\rmc}$, $A = 0$ is the unique fixed point of Eq.~\eqref{eq:OmegaA} since $\epsilon < 0$. The positive fixed point $A_{*}$ of Eq.~\eqref{eq:OmegaA} below the critical temperature $\beta \gtrsim \beta_{\rmc}$ is given by
\begin{align}
    A_{*} \simeq \sqrt{\frac{3 \epsilon}{2(r+3)}},
\end{align}
which is consistent with the critical scaling \eqref{eq:msp} of the magnetization since $m = A + 2 \rho = A + O(A^3)$. Note that
\begin{align}
    \Delta = - A - \frac{r+2}{2 r^2} \epsilon A + \frac{r+2}{r^2} A^3 + O(A^5)
\end{align}
remains bounded in time as long as the amplitude $A$ is bounded.

The amplitude $A$ relaxes exponentially toward its stationary value, $\Omega \simeq \frac{\epsilon}{r} A$ for $\beta \lesssim \beta_{\rmc}$ and $\Omega \simeq \frac{2 \epsilon}{r} (A_* - A)$ for $\beta \gtrsim \beta_{\rmc}$. The relaxation times $\tau_{+}$ and $\tau_{-}$ above and below the critical temperature show the critical slowing down with the identical exponent,
\begin{align}
    \tau_{+} = 2 \tau_{-} = \frac{r}{|\epsilon|} \sim |\beta - \beta_{\rmc}|^{-1}.
\end{align}
At the critical point $\epsilon = 0$, the relaxation time diverges and the magnetization exhibits the algebraic decay,
\begin{align}
    |m(t)| \simeq |A(t)| \simeq \sqrt{\frac{3r}{4(r+3)t}} \sim t^{- \frac{1}{2}}
\end{align}
for $t \gg 1$.

\section{Conclusion}
\label{sec:conclusion}

We study the Gibbsian structure and the detailed balance condition in stationary Nash equilibria of finite-state mean-field games with entropic control costs. We introduce a mean-field Ising game in which each player controls the flip rate of an Ising spin state to balance a ferromagnetic interaction cost against a relative-entropy control cost. The game is described by a two-state MFG that can be reduced to the two-dimensional dynamical system.

At stationary Nash equilibria, ordered equilibria bifurcate continuously from the disordered equilibrium at $\beta_{\rm c}=2(r+2)/J$, with $m\sim(\beta-\beta_{\rm c})^{1/2}$. The analysis of non-stationary MFG Nash equilibria near the critical temperature reveals the critical scaling of the relaxation time $\tau \sim |\beta - \beta_{\rmc}|^{-1}$.

In conclusion, we offer a brief comment on future directions for this study. The finite-size effect in the finite-player game is an interesting problem. In the mean-field Ising model, the magnetization distribution becomes bimodal below the critical temperature, and the stochastic spin dynamics exhibits rare transitions between metastable states. A detailed analysis of the finite-player game could clarify the role of finite population fluctuations, metastability, and equilibrium selection in the mean-field Ising game.

\begin{acknowledgements}
	This work is funded by Toyota Konpon Research Institute, Inc., JSPS KAKENHI [grant nos. 25K17309 (K.H.), 26K17606 (H.K)]
\end{acknowledgements}

\appendix

\section{On entopic cost}
\label{app:entcost}

We show that the entropic running cost function \eqref{eq:entcost} is the instantaneous relative entropy rate of an individual player's controlled process with respect to the uncontrolled process conditioned on the initial state. The transition probability density function induced by the controlled rate $\lambda$ is
\begin{widetext}
\begin{align}
    p_{\bm{\alpha}}(t+dt,\bms^{\prime}|t,\bms) =
    \begin{cases}
        \lambda_i (t,s,s^{\prime}) dt + o(dt) & \text{if $s_i \neq s_i^{\prime}$ and $\bms_{-i} = \bms^{\prime}_{-i}$ for some $i \in [N]$}, \\
        1 - \sum_{j \in [N]} \sum_{\sigma \in E \backslash \{s_j\}} \lambda_j (t,s_j,\sigma) dt + o(dt) & \text{if $\bms = \bms^{\prime}$}, \\
        o(dt) & \text{otherwise}
    \end{cases}
\end{align}
Hence, the instantaneous relative entropy rate of $p_{(\bm{\alpha}_{-i},\alpha)}(t+dt,\cdot|t,\bms)$ with respect to $p_{(\bm{\alpha}_{-i},0)}(t+dt,\cdot|t,\bms)$ is given by
\begin{align}
    &\lim_{h \downarrow 0}\frac{D(p_{(\bm{\alpha}_{-i},\alpha)}(t+h,\cdot|t,\bms)\|p_{(\bm{\alpha}_{-i},0)}(t+h,\cdot|t,\bms))}{h}
    \notag \\
    &\coloneqq \lim_{h \downarrow 0} \frac{1}{h} \sum_{\bms^{\prime} \in E^N} p_{(\bm{\alpha}_{-i},\alpha)}(t+h,\bms^{\prime}|t,\bms) \ln \frac{p_{(\bm{\alpha}_{-i},\alpha)}(t+h,\bms^{\prime}|t,\bms)}{p_{(\bm{\alpha}_{-i},0)}(t+h,\bms^{\prime}|t,\bms)}
    \notag \\
    &= \sum_{s^{\prime} \in E \backslash \{s_i\}} \left( \lambda_i (t,s_i,s^{\prime}) \ln \frac{\lambda_i (t,s_i,s^{\prime})}{k(s_i,s^{\prime})} - \lambda_i (t,s_i,s^{\prime}) + k(s_i,s^{\prime}) \right)
\end{align}
With the multiplicative parametrization \eqref{eq:lambda}, it becomes
\begin{align}
    \sum_{s^{\prime} \in E \backslash \{s_i\}} \left( e^{\alpha(t,s_i,s^{\prime})} \alpha (t,s_i,s^{\prime}) - e^{\alpha (t,s_i,s^{\prime})} + 1 \right) k(s_i,s^{\prime}).
\end{align}
This shows Eq.~\eqref{eq:entcost} is equal to the relative entropy rate of the controlled process with respect to the uncontrolled process with $k(x,y) = 1$.

\section{Derivation of coupled HJB equation}
\label{app:cHJB}

We derive the coupled HJB equation used in the main text. Fix the Markov feedback controls of all players except the $i$-th player, denoted by $\bm\alpha_{-i}$. For a test function $u:E^N\to\mathbb{R}$, we define the Bellman Hamiltonian of the $i$-th player,
\begin{align}
    G_i(t,\bms,u,\alpha_i;\bm\alpha_{-i}) &\coloneqq c_{i1}(\bms)
    +2\beta^{-1}\sum_{s^{\prime}\in E\backslash\{s_i\}}
    \left(e^{\alpha_i(s_i,s^{\prime})}\alpha_i(s_i,s^{\prime})
    -e^{\alpha_i(s_i,s^{\prime})}+1\right)
    \notag\\
    &+\sum_{s^{\prime}\in E\backslash\{s_i\}}
    e^{\alpha_i(s_i,s^{\prime})}
    \left[u(\bms_{-i},s^{\prime})-u(\bms)\right]
    +\sum_{j\in[N]\backslash\{i\}}\sum_{s^{\prime}\in E\backslash\{s_j\}}
    e^{\alpha_j(t,s_j,s^{\prime})}
    \left[u(\bms_{-j},s^{\prime})-u(\bms)\right].
\end{align}
For the value function $v_i^N$, the dynamic programming principle gives a HJB equation,
\begin{align}
    \partial_t v_i^N-rv_i^N+
    \min_{\alpha_i}G_i(t,\bms,v_i^N,\alpha_i;\bm\alpha_{-i})=0.
\end{align}
The minimization over $\alpha_i$ is pointwise in each possible jump. We notice that, for any $\alpha,\delta\in\mathbb{R}$,
\begin{align}
    e^\alpha\alpha-e^\alpha+1+e^\alpha\delta\geq 1-e^{-\delta},
\end{align}
with equality at $\alpha=-\delta$. Taking
\begin{align}
    \delta=\frac{\beta}{2}\left[u(\bms_{-i},s^{\prime})-u(\bms)\right],
\end{align}
we obtain
\begin{align}
    \min_{\alpha_i} G_i(t,\bms,u,\alpha_i ; \bm{\alpha}_{-i}) &= c_{i1}(\bms) + 2 \beta^{-1} \sum_{s^{\prime} \in E \backslash \{s_i\}} \left( 1 - e^{- \frac{\beta}{2} \left( u((\bms_{-i},s^{\prime})) - u(\bms) \right)} \right)
    \notag \\
    &\phantom{=} +\sum_{j\in[N]\backslash\{i\}}\sum_{s^{\prime}\in E\backslash\{s_j\}}
    e^{\alpha_j(t,s_j,s^{\prime})}
    \left[u(\bms_{-j},s^{\prime})-u(\bms)\right].
\end{align}
The minimizer of $G_i$ is given by
\begin{align} \label{eq:minimizer}
    \alpha_i^*(s_i,s^{\prime})
    =-\frac{\beta}{2}\left[u(\bms_{-i},s^{\prime})-u(\bms)\right].
\end{align}
A verification argument in optimal control theory ensures that the best response control is given by Eq.~\eqref{eq:minimizer} with $u = v_i^N$, which yields \eqref{eq:finiteNash}. Finally, imposing the fixed-point condition that every other player $j$ uses the corresponding minimizing feedback control gives Eq.~\eqref{eq:cHJB}.

\section{Heuristic derivation of mean-field game}
\label{app:MFG}

We provide a heuristic derivation of the MFG system. For large $N$, we suppose the functional ansatz
\begin{align} \label{eq:functionalansatz}
    v_i^N(t,\bms) = U(t,s_i,\rho^N_{\bms,i}) + O(N^{-1}).
\end{align}
A direct calculation gives
\begin{align}
    \rho^N_{(\bms_{-j},s^{\prime}),i}(s)-\rho^N_{\bms,i}(s)
    =\frac{1}{N-1}\left(\delta_{s,s^{\prime}}-\delta_{s,s_j}\right).
\end{align}
Hence,
\begin{align}
    &v_i^N(t,(\bms_{-j},s^{\prime}))-v_i^N(t,\bms)
    \notag\\
    &\quad=\frac{1}{N-1}\sum_{s^{\prime\prime}\in E}
    \left(\delta_{s^{\prime\prime},s^{\prime}}-\delta_{s^{\prime\prime},s_j}\right)
    \frac{\delta U}{\delta\rho(s^{\prime\prime})}
    (t,s_i,\rho^N_{\bms,i})+O(N^{-2}).
\end{align}
Substitution into the coupled HJB equation and replacement of sums over the other players by their empirical distribution give the equation for $U$,
\begin{align} \label{eq:master}
    0={}&(\partial_t-r)U(t,s,\rho) + \bar{c}(s,\rho)
    +2\beta^{-1}\sum_{s^{\prime} \in E\backslash\{s\}}
    \left[1-e^{-\frac{\beta}{2}(U(t,s^{\prime},\rho)-U(t,s,\rho))}\right]
    \notag\\
    &+\sum_{s^{\prime},s^{\prime \prime}\in E}
    \left[\rho(s^{\prime})e^{-\frac{\beta}{2}(U(t,s^{\prime \prime},\rho)-U(t,s^{\prime},\rho))}
    -\rho(s^{\prime \prime})e^{-\frac{\beta}{2}(U(t,s^{\prime},\rho)-U(t,s^{\prime \prime},\rho))}\right]
    \frac{\delta U}{\delta\rho(s^{\prime \prime})}(t,s,\rho).
\end{align}
We remark that $\rho$ is merely a parameter in the above equation.

We can obtain an approximate Nash equilibrium in the finite-player game by reducing the problem to the characteristics of $U$, without directly solving the above equation for $U$. We now suppose that the initial values $\bmS (0)$ are sampled independently according to the one-body probability density function $\rho_0$. The forward master equation for the $N$-body probability density function $\rho^N(t,\bms)$ under a Nash equilibrium of the finite-player game is
\begin{align}
    \partial_t\rho^N(t,\bms)
    = \sum_{i=1}^N\sum_{s^{\prime}\in E\backslash\{s_i\}}
    \left[\rho^N(t,(\bms_{-i},s^{\prime}))
    e^{-\frac{\beta}{2}\left(v_i^N(t,\bms)-v_i^N(t,(\bms_{-i},s^{\prime}))\right)} -\rho^N(t,\bms)
    e^{-\frac{\beta}{2}\left(v_i^N(t,(\bms_{-i},s^{\prime}))-v_i^N(t,\bms)\right)}\right]
\end{align}
\end{widetext}
with the initial data $\rho^N(0,\bms) = \prod_{i=1}^{N} \rho_0 (s_i)$. Summing the above differential equation over $s_2,\dots,s_N$ provides the equation for the one-body marginal density function $\rho_1^N(\bms_1) = \sum_{s_2,\dots,s_N} \rho^N(\bms)$. Assuming the conditional law of large numbers for the empirical density function, for a continuous function $c$ of the empirical density function, we can evaluate
\begin{align}
    &\sum_{s_2,\ldots,s_N}\rho^N(t,\bms) c(s_1,\rho^N_{\bms,1})
    \notag \\
    &= \rho_1^N(t,s_1) \bbE \left[ c(s_1,\rho^N_{\bmS,1}) | S_1 = s_1 \right]
    \notag \\
    &\approx \rho (t,s_1) c(s_1,\rho (t)),
\end{align}
where $\rho \coloneqq \lim_{N \uparrow \infty} \rho^N_1$. Applying this replacement to the one-body marginal of the $N$-body master equation gives Eq.~\eqref{eq:MFmaster} with $U(t,s,\rho(t))$ in place of $v(t,s)$.

We now define
\begin{align}
    v(t,s)\coloneqq U(t,s,\rho(t)).
\end{align}
The chain rule gives
\begin{align}
    \partial_t v(t,s)
    =\partial_t U(t,s,\rho(t))
    +\sum_{s^{\prime}\in E}\partial_t\rho(t,s^{\prime})
    \frac{\delta U}{\delta\rho(s^{\prime})}(t,s,\rho(t)).
\end{align}
The second term exactly cancels the term in the second line of Eq.~\eqref{eq:master} when Eq.~\eqref{eq:MFmaster} with $U(t,s,\rho(t))$ in place of $v(t,s)$ is used. This yields the HJB equation \eqref{eq:MFHJB}, and hence the MFG system stated in the main text.

Finally, when the initial values $\bmS (0)$ are sampled independently according to $\rho_0$, the Nash equilibrium of finite-player game along typical trajectories is then approximated in the limit $N \uparrow \infty$ by
\begin{align}
    &\varphi_i^{N,*}(t,\bmS(t))(S_i(t),s^{\prime})
    \notag \\
    &\approx -\frac{\beta}{2}\left[U(t,s^{\prime},\rho^N_{\bmS(t),i})
    -U(t,S_i(t),\rho^N_{\bmS(t),i})\right]
    \notag\\
    &\approx -\frac{\beta}{2}\left[U(t,s^{\prime},\rho(t))
    -U(t,S_i(t),\rho(t))\right]
    \notag\\
    &=-\frac{\beta}{2}\left[v(t,s^{\prime})-v(t,S_i(t))\right],
\end{align}
In the second line, we use the functional ansatz. In the third line, we use the law of large numbers for the empirical density function, $\rho^N_{\bmS(t),i} \to \rho (t)$ with high probabilities. In the fourth line, we use the definition of $v$.

\section{Boundedness of value function}
\label{app:bound}

We show that the value function is bounded in time. Suppose that $\bar{c}$ is a continuous function of $\rho$. Then, there exists a constant $C_0$ such that $|\bar{c}(s,\rho)| \leq C_0$ for any $s \in E$ and $\rho$ because the set of probability density functions on $E$ is compact. Since $c_2(s,\alpha) \geq 0$,
\begin{align}
    C(t,s,\alpha) \geq \int_t^{\infty} du \ e^{-r(u-t)} (-C_0) = - \frac{C_0}{r}
\end{align}
for all controls $\alpha$. Hence, $v(t,s) \geq - C_0/r$. On the other hand,
\begin{align}
    v(t,s) \leq C(t,s,\alpha = 0) \leq \int_t^{\infty} du \ e^{-r(u-t)} C_0 = \frac{C_0}{r}.
\end{align}
Therefore, $|v(t,s)| \leq C_0/r$.

\section{Solutions of Eq.~\eqref{eq:Fm}}
\label{app:beta}

We define
\begin{align}
    \hat{\beta}(m) \coloneqq \frac{2}{J m} \left( \frac{r}{2} \ln \frac{1+m}{1-m} + \frac{2m}{\sqrt{1-m^2}} \right).
\end{align}
The positive solutions of Eq.~\eqref{eq:Fm} are given by the solutions of $\beta = \hat{\beta}(m)$ for $m \in (0,1)$. We show that $\hat{\beta}(m)$ is strictly increasing in $(0,1)$. By differentiating $\hat{\beta}$, we obtain
\begin{align}
    \hat{\beta}^{\prime}(m) = \frac{2}{Jm} F_{\hat{\beta}(m)}^{\prime}(m).
\end{align}
A similar argument below Eq.~\eqref{eq:Fprime} implies $\hat{\beta}^{\prime}(m) > 0$ for $m \in (0,1)$. Notice that
\begin{align}
    \lim_{m \downarrow 0} \hat{\beta}(m) = \beta_{\rmc}, \quad \lim_{m \uparrow 1} \hat{\beta}(m) = \infty.
\end{align}
Hence, a positive solution exists if and only if $\beta > \beta_{\rmc}$, and the positive solution $m_{*}(\beta)$ is unique if it exists. Since $F_{\beta}$ is an odd function, $F_{\beta}(m) = 0$ has the unique solution $m=0$ for $\beta \leq \beta_{\rmc}$ and has exactly three solutions, $m = 0$ and $m = \pm m_{*}(\beta)$ for $\beta > \beta_{\rmc}$.

\bibliographystyle{jabbrv_apsrev4-2}
\bibliography{myref}

\end{document}